\documentclass[12pt]{article}
\usepackage{geometry}                % See geometry.pdf to learn the layout options. There are lots.
\usepackage{graphicx}
\usepackage{amssymb}					
\usepackage{amsmath}					
\usepackage{mathrsfs}
\usepackage{tikz}
\usepackage{epstopdf}
\usepackage{amsthm}
\newtheorem*{theorem*}{Theorem}
\newtheorem*{lemma*}{Lemma}
\newtheorem*{example*}{Example}
\newtheorem{theorem}{Theorem}

\def\b0{{\bf 0}}
\def\b1{{\bf 1}}

\def\cC{{\cal C}}

\def\cC{{\cal C}}

\def\n{\noindent}

\begin{document}
\title{
%Which embeddings have Hamiltonian extensions?
Graph embeddings with no Hamiltonian extensions \thanks{to appear in {\it Bull. of the Inst. of Comb. \& its Appl.}, accepted for publ. 2/27/2023}}
\author{Paul C. Kainen $\;\;\;\;\;\;\;\;\;\;\;\;$ Shannon Overbay\\ \texttt{kainen@georgetown.edu $\;\;\;\;\;$overbay@gonzaga.edu} }
\date{}                                           % Activate to display a given date or no date

\newcommand{\Addresses}{{% additional braces for segregating \footnotesize
  \bigskip
  \footnotesize

%\n Paul C. Kainen, \textsc{Department of Mathematics and Statistics,\\ Georgetown University, 37th and O Streets, N.W., Washington DC 20057}\\
\vspace{-0.27cm}

\n

\par\nopagebreak
}}

\maketitle
%\section{}
%\subsection{}

\abstract{\n 
{\it We show that extending an embedding of a graph $\Gamma$ 
in a surface to an embedding of a Hamiltonian  supergraph can be blocked by 
%the existence of 
certain planar subgraphs but, for some subdivisions of $\Gamma$, Hamiltonian extensions must exist.
%there are Hamiltonian extensions of a subdivision of $\Gamma$.
}  
%{\it A {\rm (2-cell)} embedding of a graph in a surface is {\rm extended} by including new vertices and non-crossing edges in some of the open-disk regions; extension is {\rm Hamiltonian} if the resulting supergraph is.  A family of embeddings {\rm (obstacles)} which do not have Hamiltonian extensions is shown, and we conjecture that these obstacles are the {\rm only} non-Hamiltonian-extendable embeddings. In contrast, we prove that there are no obstacles for the topological variant of the problem which allows one to take an arbitrary subdivision of the edges of the original embedded graph before trying to extend the embedding.}
}
\smallskip

\n
{\bf Key Phrases}: {\it extending embeddings, Hamiltonian cycle in embedded graph}.

%\subsubsection*
%{A geometric proof that $mgt(G \,\Box\, H) \leq mgt(G) + mgt(H)$ }
\bigskip

\section{Introduction}

The objects studied in this paper are 2-cell embeddings of graphs in (closed) surfaces.  We ask: {\it When can such an embedding be extended to an embedding of a Hamiltonian graph, containing the original graph as a subgraph?} The embedding is into the same surface so that the supergraph is obtained as a subdivision of some of the regions of the original embedding but the edges of the original graph are not subdivided.  See Fig. 1 below.

This problem is a variant of the differently specified
%more specialized 
question asked in \cite{mo}, ``{\it When
is a graph, embeddable on a surface S, a subgraph of a Hamiltonian graph
which is also embeddable on S?}''  McKenzie and Overbay showed \cite{mo} that the bipartite complete graphs, with genus $\gamma \leq 1$ which are {\it not} Hamiltonian, are subgraphs of genus-$\gamma$  graphs that {\it are} Hamiltonian.  

The formulation here emphasizes the  embedding itself, rather than the possibility of being embedded.
The idea of extending graph invariants to graph embeddings goes back (at least) to \cite{cook-1, cook-2, k73, k-abh, pz1988}.

Merely being non-Hamiltonian isn't enough to prevent a Hamiltonian extension.  For instance, the Petersen graph has an embedding in the torus, and one can add three edges to the embedding to make the enlarged graph Hamiltonian where each added edge occurs within a region of the original embedding.
{\it Which embeddings ensure that no such Hamiltonian extension can be found?}

We obtain a large family of non-Hamiltonian-extendable embeddings using an idea of Klee (see Malkevitch \cite{malkevitch}) and conjecture that there are no other such non-Hamiltonian-extendable embeddings.  However,
if the edges of the original graph can be subdivided before trying to extend it, then we show that every graph embedding has such a {\it topological} Hamiltonian extension.  

The paper proceeds as follows: Section 2 has definitions; in Section 3 we build non-Hamiltonian-extendable graph embeddings.  Section 4 proves that weakening the condition of extendability to allow subdivision of edges of the original graph makes it possible to always find a Hamiltonian extension.

\section{Definitions}
%Recall that a {\bf stellation} of a triangle in a planar surface adds a new vertex in the center with edges to the three corners. {\bf Twice stellate} means apply stellation to the stellation.

A {\bf 2-cell embedding} $i$ of a finite graph $\Gamma$ in a surface $S$ is a continuous embedding $i: \Gamma \to S$ such that $S \setminus i(\Gamma)$ is a disjoint union of open 2-disks, the {\bf regions} (of $i$).
%; an embedding is a {\bf closed 2-cell embedding} if the closure of each open region is a closed 2-disk whose boundary is a cycle contained in the graph.  
%If $\Gamma'$ is a subdivision of $\Gamma$, then $i$ induces a corresponding 2-cell embedding $i'$ of $\Gamma'$ in $S$.  
If $G$ is some graph which contains $\Gamma$ as a subgraph and 
$j: G \to S$ is a 2-cell embedding, then we say that {\bf j extends i} if $j |_\Gamma = i$.  
%If $i$ is a closed 2-cell embedding, we require this of $j$ as well.

We call a 2-cell embedding $i$ of $\Gamma$ in $S$ {\bf Hamiltonian extendable} if $i$ can be extended to an embedding of a Hamiltonian supergraph $G$ in $S$.  
Otherwise, $i$ is {\bf non-Hamiltonian-extendable}.

\begin{center}
\includegraphics[width=0.7\linewidth]{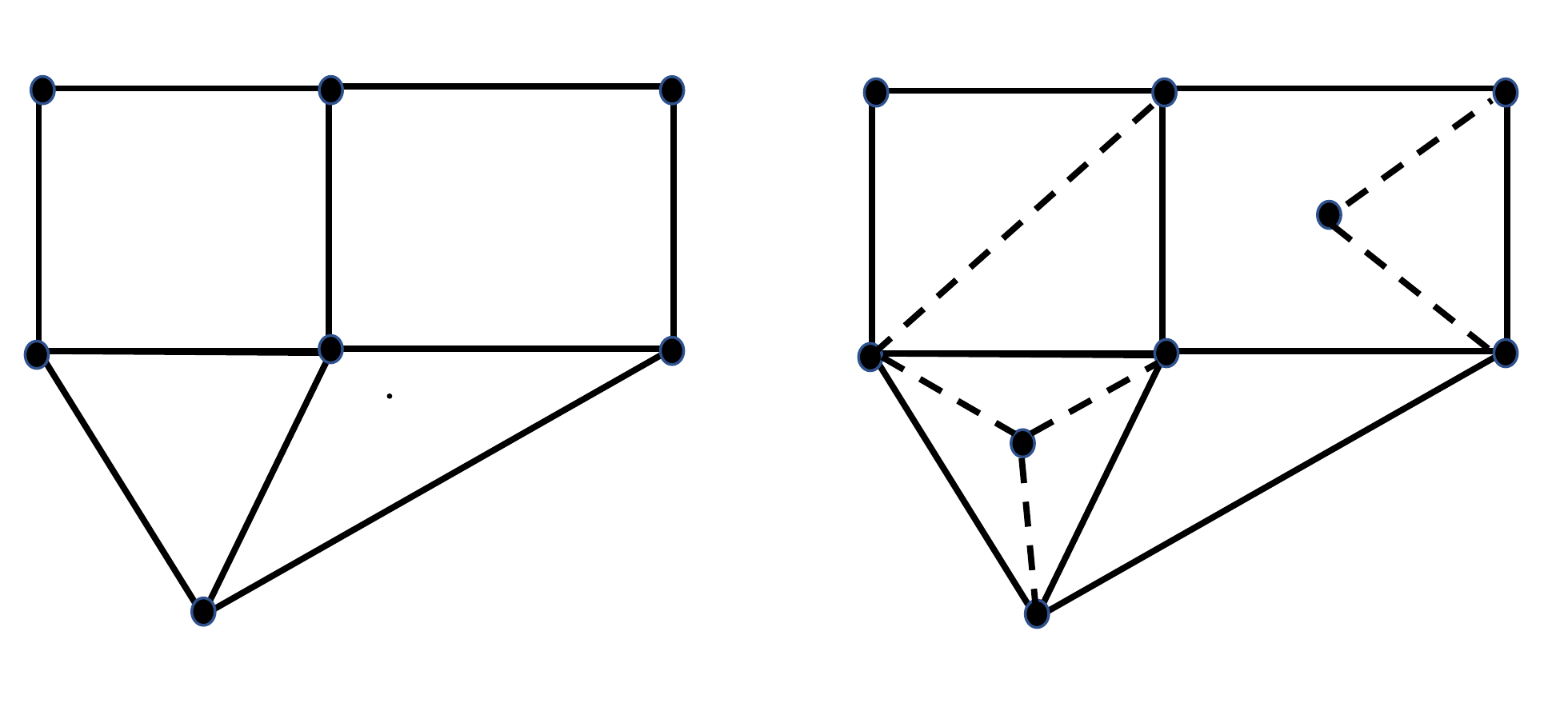}
\end{center}
\centerline{{\bf Figure 1} \ \ {\it An extension of an embedding}}
\label{fig:1st}
\bigskip
\bigskip

A path or cycle is {\bf oriented} if its edges are assigned a consistent direction.  If $P$ is an oriented path, let $P^o$ denote $P$ {\it minus its terminal point}. 
A subdivision of an edge is a path whose endpoints agree with the endpoints of the edge.  A {\bf subdivsion} of a graph is a graph obtained by subdividing some or all of the edges.  Two graphs are homeomorphic iff they have isomorphic subdivisions.

An embedding $i:\Gamma \to S$ will be said to have a {\bf topological extension} if there exists a subdivision $\Gamma'$ of $\Gamma$ and an extension $j: G \to S$ of $i'$, where $i'$ is the embedding $\Gamma' \to S$ induced by $i$.

A 2-cell embedding $i: \Gamma \to S$ is of {\bf Klee type} if the number $r$ of regions exceeds the number $p$ of vertices; $i$ is of {\bf local} Klee type if there exists a cycle $C$ contained in $\Gamma$ such that (i) $C$ separates $\Gamma$, (ii) $i(C)$ separates $S$ (into {\it inside} and {\it outside}), and (iii) 
if $r_C$ is the number of regions of $i(\Gamma)$ inside $C$ and $p_C$ is the number of vertices of $\Gamma$ inside or on $C$, then $r_C \geq p_c$. See Fig. 2. Labeling  inside/outside is arbitrary and both parts of $S \setminus C$ could be nonplanar.

\section{Graph embeddings of Klee type }

Extending a 2-cell embedding of Klee or local Klee type to include points in the interiors of too many regions must produce a non-Hamiltonian-extendable graph.
We conjecture that these obstacles are the only way to produce such non-Hamiltonian-extendable graphs.

\begin{theorem}
(a) Let $i:\Gamma \to S$ be an embedding of Klee type with $r > p$. Then, for any extension $j:G \to S$, $G$ is not Hamiltonian provided $G$ contains vertices $w_1, \ldots, w_s$ inside distinct regions of $i$, $R_1, \ldots, R_s$, for $r \geq s \geq p+1$.\\
(b) Let $i:\Gamma \to S$ be an embedding of local Klee type with $r_C \geq p_C$. Then, for any extension $j:G \to S$, $G$ is not Hamiltonian provided $G$ contains vertices $w_1, \ldots, w_s$ inside distinct regions of $i$, $R_1, \ldots, R_s$, inside $i(C)$ for $r_C \geq s \geq p_C$.
%Then any embedding $i':\Gamma' \to S$ obtained from $i$ by stellating $s > p$ has no Any graph embedding which extends the stellation of a graph embedding with the Klee property is not Hamiltonian.
\end{theorem}
\begin{proof}
%Observe that to have a cycle through $w_k$ it must have degree at least 2 in $G$.
We argue by contradiction. Suppose there is an extension $j: G \to S$ of $i$ and let
$Z$ be any oriented cycle contained in $G$ which includes all $s$ points.  By construction, between any two consecutive (with respect to $Z$) points, say, $w_k, w_{k+1}$ ($k = 1, \dots, s$, addition mod $s$), there is a unique vertex $v_k$ in the boundary of the region $R_k$ of $i$ containing $w_k$ such that $v_k$ is in $Z$ and the subpath $P_k$ of $Z$ from $w_k$ to $v_k$ contains no other point in $W:=\{w_1,\ldots,w_s\}$ and no other point in the boundary of $R_k$.   

In case (a), $Z$ contains at least $s$ points in $V\Gamma$, which contradicts the assumption $s \geq p+1$. In case (b), if $r_C > p_C$, then as in (a), no such cycle $Z$ can exist, while if $r_C=p_C$, the only possibility is that $Z$ includes all vertices on $C$ (and some inside it), so $Z$ can't include the vertices of $\Gamma$ outside $C$.
\end{proof}
%A stellation of a 2-cell graph embedding is any embedding obtained by introducing a vertex in the interior of one or more of the regions along with edges joining this vertex to all the vertices in the bounding walk of that region. A stellation is total if each region gets such an interior cone-point.

%As an example of a non-Hamiltonian-extendable graph, one may take the sphere and form its total stellation, which has 5 vertices and 6 regions. This is a graph embedding of Klee type, so its total stellation is not Hamiltonian.

Using the genus formula \cite{ringel, beineke-harary} for cubes, $\gamma(Q_d) = 1+(d-4)2^{d-3}$, easy calculation shows that for the $d$-cube, the number of regions in the genus embedding  is $r := r(d) := d \,2^{d-2} > 2^d = p$ for $d \geq 5$.  Indeed, by Euler's formula, 
%we have the equation
%\begin{equation}
\[
2^d -d2^{d-1} + r(d) = 2 - 2(d-4)2^{d-3} - 2.
\]
%\end{equation}
Solving for $r(d)$ gives the result.
So cubes of dimension $\geq 5$ are of Klee type.  Using the construction in case (a) above, one obtains for the 5-cube, by adding one new vertex $w$ in the middle of $s$ of the square faces, $33  \leq s \leq 40$, and using any of the 11 ways to connect each $w_k$ to $\geq 2$ of the 4 vertices on the boundary of the face which contains it, the number of distinct 2-connected non-Hamiltonian graphs with embedding in $S_5$ extending that of the 5-cube is
$$N = \sum_{k=33}^{40} {{40}\choose{k}} 11^k \approx 1.45 \times 10^{43}.$$
%Similarly, $K_n$ is of Klee type, for $n \geq 7$ and $n \equiv 0,3,4,7$ (mod 12). 
%The hyperoctahedra are left as an exercise for the reader.

The {\bf stellation} of a triangular region puts one new vertex into the interior and joins it to all three corners.
Iterating this operation on the resulting three triangles gives a local Klee type graph embedding with $C = K_3$, where $r_C=9$ and $p_C=7$.  Hence, stellating all 9 of the regions produces a non-Hamiltonian-embeddable graph, no matter where it occurs in some potentially large graph embedding.  Here the inside region is what was inside the triangle.  See Fig. 2.

\begin{center}
\includegraphics[width=0.7\linewidth]{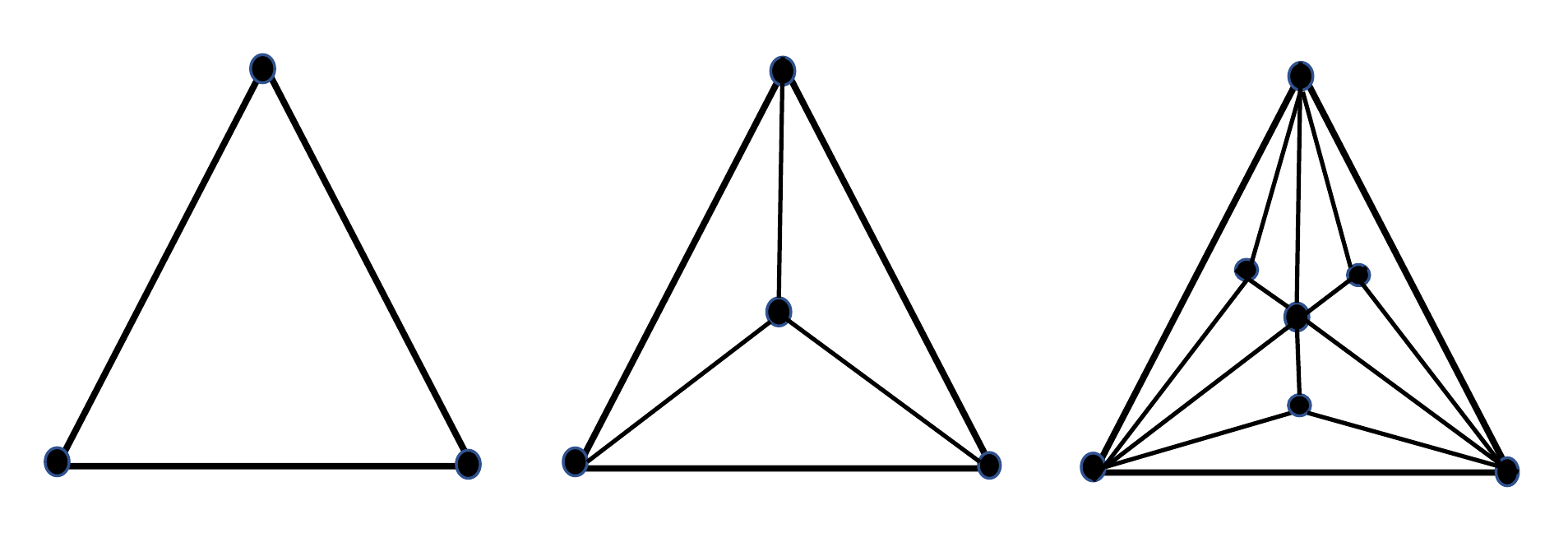}
\end{center}
\centerline{{\bf Figure 2} \ \ {\it Local Klee type embedding on right; triangle C on left}}
\label{fig:1st}
\bigskip
\bigskip

\section{Topological extensions }

The planar case of Theorem \ref{th:ham-top} below is (implicitly) in \cite[p. 32]{so-thesis}.

\begin{theorem}
Any embedding $\Gamma \subset S$ has a Hamiltonian topological extension.
\label{th:ham-top}
\end{theorem}
\begin{proof}
Let $i:\Gamma \to S$ be an embedding.  Consider the $p = |V\Gamma|$ points $i(v) \in S$ for $v \in V\Gamma$.  As $S$ is a closed surface, it cannot be disconnected by the removal of any path (or any other contractible subset).  Hence, for any enumeration of the points $i(v)$, say $i(v_1), \ldots, i(v_p)$, there is a topological path $P_1$ in $S$ from $v_1$ to $v_2$, then a path $P_2$ in $S \setminus P_1^o$ from $v_2$ to $v_3$, and so on, until one chooses a path $P_p$ in $S \setminus \bigcup_{k=1}^{p-1} P^o_k$ from $i(v_p)$ to $i(v_1)$.  
The union of the paths $P_1, \ldots, P_p$ is a non-self-intersecting closed curve $\cC$ with $i(V\Gamma) \subset \cC \subset S$.  
Since $S$ is triangulable, there is an arbitrarily small perturbation $\cC'$ of $C \setminus i(V\Gamma)$ so that 
\[
\cC' \cap i(\Gamma) = i(V\Gamma) \cup Y,
\]
where $Y$ is a finite set of points at which $\cC'$ crosses interiors of edges of $\Gamma$.
%\newpage

%we can ensure that $\cC$ intersects $i(\Gamma)$ only in isolated points which consist of the vertices of $\Gamma$ together with some points where $\cC$ crosses the interior of an edge of $\Gamma$.
Take the points in $Y$
%these latter intersection points 
as subdivision vertices for the edges of $\Gamma$, and let $\Gamma'$ be the resulting subdivision of $\Gamma$.  Define a graph $G$ as the union of $\Gamma'$ and the new edges which result by subdividing $\cC'$ using both the vertices of $\Gamma'$ and the subdivision points.  The resulting copy of $G$ in $S$ extends the embedding of $\Gamma'$ and $G$ has the subdivided $\cC'$ as Hamiltonian cycle.
\end{proof}

%\section{Remarks}
We ask: What is the least number of subdivision points needed?

An alternate means to find a Hamiltonian embedding extending a subdivision of some given embedding might be achievable using the ``mesh surface'' methods in Akleman et al. \cite{wire-sculpture}.

%We found an interesting paper showing a multiplicity of Hamiltonian graphs, densely covering surfaces, which could give an alternate means to find a Hamiltonian embedding extending the given embedding; see Akleman et al. \cite{wire-sculpture}.

%For example, one can add three edges to the embedding of the Petersen graph into the torus, producing a supergraph of the Petersen graph which is topologically embedded in the torus.

\end{document}